\documentclass[12pt]{article}

\usepackage[latin1]{inputenc}
\usepackage{amsmath,amsthm,amsfonts,amscd,amssymb,xspace,color,rotating}
\usepackage[french]{babel}

\newcommand{\td}[2]{\xrightarrow[#1\rightarrow #2]{}}

\newcommand{\abs}[1]{\left|#1\right|}

\def\card{\qopname\relax o{card}}

\newtheorem{Thm}{Th\'eor\`eme}[section]
\newtheorem{Prop}[Thm]{Proposition}
\newtheorem{Lem}[Thm]{Lemme}
\newtheorem{Cor}[Thm]{Corollaire}

\newtheorem*{Thm*}{Th\'eor\`eme}
\newtheorem*{Cor*}{Corollaire}
\newtheorem*{Prop*}{Proposition}

\newcounter{ploum}

\newcounter{ex}[section]
\newcounter{rem}[section]

\numberwithin{equation}{section}

\begin{document}

\newenvironment{rem}
{\addtocounter{rem}{1}\setlength{\topsep}{1em}\par\trivlist\item{\bf Remarque
\arabic{section}.\arabic{rem}.} }{\endtrivlist}

\newenvironment{ex}
{\addtocounter{ex}{1}\setlength{\topsep}{1em}\par\trivlist\item{\bf Exemple
\arabic{section}.\arabic{ex}.} }{\endtrivlist}

\newenvironment{demo}[1][D\'emonstration]
{\setlength{\topsep}{1em}\par\trivlist\item{\em #1.} }
{\openbox\endtrivlist}

\newenvironment{cond}{\begin{list}{{\it (\roman{ploum})}}{\usecounter{ploum}}}
{\end{list}}

\newenvironment{souscond}{\begin{list}{{\hspace{5ex}\it (\alph{ploum})}}{\usecounter{ploum}}}
{\end{list}}

\begin{center}

{\Huge\bf Percolation sur le système à trois points}

\ \\

{\Large\bf J.-F. Quint}

\ \\

\end{center}

\section{Introduction}

Le système à trois points est l'ensemble  $X$ constitué des $(x_{k,l})_{(k,l)\in\mathbb Z^2}$ dans $(\mathbb Z/2\mathbb Z)^{\mathbb Z^2}$ tels que, pour tout $(k,l)$ dans $\mathbb Z^2$, on ait $x_{k,l}+x_{k+1,l}+x_{k,l+1}=0$ (dans $\mathbb Z/2\mathbb Z$), muni de l'action naturelle de $\mathbb Z^2$ par translations sur les coordonnées. Cette action est engendrée par les deux applications commutantes $T:(x_{k,l})\rightarrow (x_{k+1,l})$ et $S:(x_{k,l})\rightarrow (x_{k,l+1})$. On munit $X$ de sa topologie naturelle d'espace compact : le groupe $\mathbb Z^2$ agit par homéomorphismes sur $X$.

On peut alors chercher à classifier les probabilités boréliennes de $X$ qui sont invariantes par l'action de $\mathbb Z^2$. C'est un problème analogue à celui, posé par Fürstenberg dans \cite{Fu}, de la classification des probabilités boréliennes du cercle qui sont à la fois invariantes par le doublement et le triplement de l'angle. Cette question a été abordée par plusieurs auteurs qui ont adapté à l'espace $X$ - et à d'autres systèmes - les méthodes partielles de résolution du problème de Fürstenberg : citons, par exemple, \cite{HMM}, \cite{Sab} et \cite{Silb}.

Dans cet article, nous proposons une nouvelle approche de ce problème, basée sur des propriétés de percolation que nous allons à présent décrire.

Notons $Y$ le fermé de $X$ constitué de l'ensemble des $(x_{k,l})$ dans $X$ avec $x_{0,0}=1$. Si $x=(x_{k,l})$ appartient à $Y$, on a $x_{1,0}+x_{0,1}=1$ et, donc, un et un seul des deux éléments $Tx$ et $Sx$ appartient à $Y$. Ceci définit une application continue $\sigma:Y\rightarrow Y$. Soit $\mu$ une probabilité borélienne sur $X$ qui est invariante par l'action de $\mathbb Z^2$ (c'est-à-dire  qu'on a $T_*\mu=S_*\mu=\mu$). Alors on a  $\mu(Y)=0$ si et seulement si $\mu$ est la masse de Dirac en la famille nulle. Si $\mu(Y)>0$, la restriction de $\mu$ à l'ensemble $Y$ est quasi-invariante par $\sigma$ (c'est-à-dire que, pour tout borélien $B$ de $Y$ avec $\mu(B)=0$, on a $\mu(\sigma^{-1}(B))=0$). Si $x$ est un point de $Y$, on définit la $\sigma$-composante de $x$ comme l'ensemble des points $y$ de la $\mathbb Z^2$-orbite de $X$ qui appartiennent à $Y$ et pour lesquels il existe des entiers naturels $p$ et $q$ avec $\sigma^p(x)=\sigma^q(y)$. En d'autres termes, la $\sigma$-composante de $x$ est l'ensemble des points $y$ de la $\mathbb Z^2$-orbite de $x$ pour lesquels les chemins tracés sur la $\mathbb Z^2$-orbite de $x$ $(\sigma^p(x))_{p\geq 0}$ et $(\sigma^q(y))_{q\geq 0}$ se rejoignent.
La $\sigma$-composante de $x$ possède une structure naturelle d'arbre : les voisins de $x$ sont $\sigma(x)$ et les éventuels éléments de $\sigma^{-1}(x)$.
On appelle $\sigma$-antécédents de $x$ les éléments de $\bigcup_{n\in\mathbb N^*}\sigma^{-n}(x)$.

Notons $Z$ l'ensemble des éléments de $Y$ qui admettent une infinité de $\sigma$-antécédents.
Le résultat principal de cet article est le

\begin{Thm*} Soit $\mu$ une probabilité borélienne sans atomes, invariante et ergodique par l'action de $\mathbb Z^2$ sur $X$.
Alors une et une seule des deux propriétés suivantes est vérifiée
\begin{cond}
\item On a $\mu(Z)=0$ et il existe un ensemble borélien $\mathbb Z^2$-invariant $B$ avec $\mu(B)=1$ tel que, pour tout $x$ dans $Y\cap B$, la $\sigma$-composante de $x$ contient tous les points de la $\mathbb Z^2$-orbite de $x$ qui appartiennent à $Y$.
\item On a $\mu(Z)>0$ et il existe un ensemble borélien $\mathbb Z^2$-invariant $B$ avec $\mu(B)=1$ ayant les propriétés suivantes :\begin{souscond}
\item pour tout $x$ dans $Y\cap B$, la $\sigma$-composante de $x$ contient des points de $Z$.
\item pour tout $x$ dans $Z\cap B$, l'ensemble $\sigma^{-1}(x)$ contient un et un seul point de
$Z$.
\item pour tout $x$ dans $Y\cap B$, la $\mathbb Z^2$-orbite de $x$ contient une infinité de $\sigma$-composantes. Plus précisément, pour toute $\sigma$-composante $C$ de la $\mathbb Z^2$-orbite de $x$,
il existe un unique entier relatif $k$ tel que $T^kS^{-k}x$ appartienne à $Z\cap C$ et l'ensemble $\{k\in\mathbb Z|T^kS^{-k}x\in Z\}$ n'est ni majoré ni minoré.
\end{souscond}\end{cond}\end{Thm*}

\begin{figure}\begin{center}\input{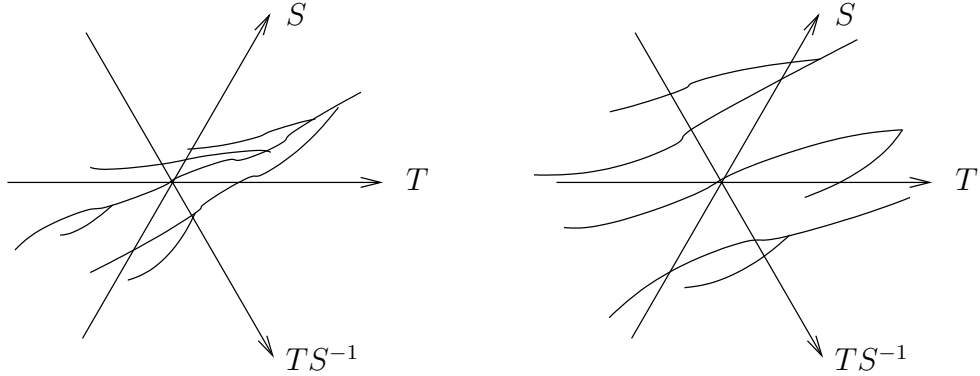}\caption{Orbites de type arbre et de type ruban}\label{arbreruban}\end{center}\end{figure}

Nous dirons qu'une mesure sans atomes $\mathbb Z^2$-invariante et
ergodique sur $X$ est du type arbre si elle vérifie la première
propriété du théorème et qu'elle est du type ruban si elle en
vérifie la seconde propriété. 
Ces deux situations sont représentées par la figure \ref{arbreruban}. Dans cette figure, à gauche, on voit, à grande distance, une orbite générique de type arbre : les points n'ont qu'un nombre fini de $\sigma$-antécédents et leurs trajectoires suivant $\sigma$ se rejoignent. \`A droite, on voit une orbite générique de type ruban : les points ayant une infinité de $\sigma$-antécédents forment des courbes transversales aux droites de la forme $\{T^kS^{-k}x|k\in\mathbb Z\}$ et les trajectoires suivant $\sigma$ des points n'ayant qu'un nombre fini de $\sigma$-antécédents rejoignent ces courbes.

\`A la section \ref{exarbre}, nous
montrerons que la mesure de Haar de $X$, vu comme sous-groupe fermé
du groupe compact $(\mathbb Z/2\mathbb Z)^{\mathbb Z^2}$, est du
type arbre. \`A la section \ref{exruban}, nous donnerons un exemple de
mesure du type ruban.

Auparavant, dans les sections \ref{antecedents} et \ref{compo}, nous
établirons un certain nombre de résultats préliminaires sur la
topologie des $\sigma$-composantes dans les $\mathbb Z^2$-orbites.
Nous conclurons la démonstration du théorème à la section
\ref{demothm}. Les résultats de la section \ref{antecedents}
reposent sur un argument analogue à celui employé par Burton et
Keane dans \cite{BK} pour établir l'unicité de la composante infinie
pour la percolation indépendante dans $\mathbb Z^d$, avec $d\geq 2$.

\section{Les $\sigma$-antécédents}
\label{antecedents}

Rappelons qu'on a noté $Z$ l'ensemble des éléments de $Y$ qui admettent une infinité de $\sigma$-antécédents.
On désigne par $Z_2$ l'ensemble des éléments $x$ de $Z$ pour lesquels l'ensemble $\sigma^{-1}(x)$ contient deux éléments
qui appartiennent tous les deux à $Z$. Dans cette section nous allons démontrer la

\begin{Prop} \label{uniqueantecedent}
Soit $\mu$ une probabilité borélienne sans atomes, invariante et ergodique par l'action de $\mathbb Z^2$ sur $X$. On a $\mu(Z_2)=0$.
\end{Prop}

Introduisons quelques notations. Soient $x$ un point de $X$ et $n$ un entier naturel. On pose
\begin{align*}\Delta_n(x)&=\{T^kS^lx|k\geq 0,l\geq 0,k+l\leq n\}\\
\Delta_n^{\circ}(x)&=\{T^kS^lx|k>0,l>0,k+l\leq n\}\\
\partial \Delta_n(x)&=\{T^kx|0\leq k\leq n\}\cup\{S^lx|0\leq l\leq n\}\cup\{T^kS^{n-k}x|0\leq k\leq n\}.\end{align*}

L'ensemble $\Delta_n(x)$ est donc un triangle de côté $n$ issu de $x$ et $\partial \Delta _n(x)$ est son bord.

Si $y$ est un point de $Y\cap \Delta_n(x)$, on définit sa
$(\sigma,x,n)$-composante comme l'ensemble des éléments $z$ de
$Y\cap \Delta_n(x)$ tels qu'il existe des entiers naturels $q$ et
$r$ avec $\sigma^q(y)=\sigma^r(z)\in \Delta_n(x)$. Par construction,
chaque $(\sigma,x,n)$-composante contient un unique point de
$\{T^kS^{n-k}x|0\leq k\leq n\}$. L'intérêt de ces définitions
provient du lemme suivant, qui est un analogue du lemme de Burton et
Keane dans \cite{BK} :

\begin{Lem} Soient $x$ un point de $X$ et $n$ un entier naturel. Si $C\subset \Delta_n(x)$ est une  $(\sigma,x,n)$-composante qui contient des points de $Z\cap \Delta_n^{\circ}(x)$,
on a $$\card(Z_2\cap C\cap \Delta_n^{\circ}(x))\leq \card(\partial \Delta_n(x)\cap C)-2.$$\end{Lem}

L'idée de la démonstration de ce lemme est que $C$ est muni d'une structure naturelle d'arbre dans laquelle les points ont deux ou trois voisins. Un point $z$  de $Z_2\cap C\cap \Delta_n^{\circ}(x)$ a trois voisins et les composantes connexes de $C-\{z\}$ rencontrent toutes $\partial \Delta_n(x)$.

\begin{figure}\begin{center}\input{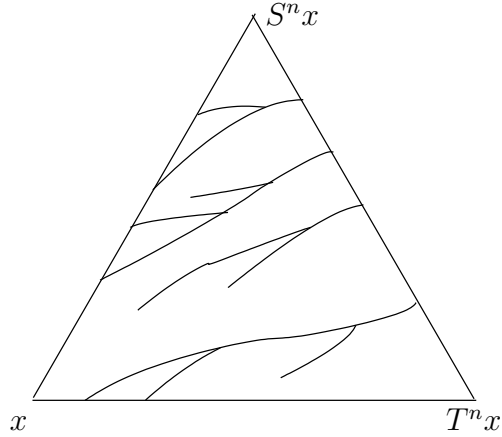}\caption{$(\sigma,x,n)$-composantes}\label{petitarbre}\end{center}\end{figure}

\begin{demo} Comme $C$ contient des points de $\Delta_n^{\circ}(x)$, $C$ contient un unique point de $\{T^kS^{n-k}x|1\leq k\leq n-1\}$. Par ailleurs, si $y$ est un point de
$C\cap Z$, l'un au moins des $\sigma$-antécédents de $y$ appartient à $\partial \Delta_n(x)$. Comme ce $\sigma$-
antécédent n'appartient pas à  $\{T^kS^{n-k}x|1\leq k\leq n-1\}$,
$C\cap \partial \Delta_n(x)$ contient au moins deux points, et le résultat est établi si $Z_2\cap C\cap \Delta_n^{\circ}(x)=\emptyset$.

Dans le cas contraire, si $y$ est un point de $Z_2\cap C\cap
\Delta_n^{\circ}(x)$, au moins deux des $\sigma$-antécédents de $y$
appartiennent à $\partial \Delta_n(x)$ et, donc, $C\cap \partial
\Delta_n(x)$ contient au moins trois points. L'ensemble $C$ possède
une structure naturelle d'arbre pour laquelle les voisins d'un
élément $z$ de $C$ sont $\sigma(z)$, si ce dernier appartient à
$\Delta_n(x)$, et les éventuels éléments de $\sigma^{-1}(z)$ qui
appartiennent à  $\Delta_n(x)$. Alors, si $z$ est un point de
$Z_2\cap C\cap \Delta_n^{\circ}(x)$, $C-\{z\}$ contient trois
composantes connexes pour la structure d'abre et ces trois
composantes connexes déterminent une partition $P_z$ de  $\partial
\Delta_n(x)\cap C$ en trois ensembles non vides. L'appplication
$z\mapsto P_z$ est injective et l'ensemble $\mathcal P=\{P_z|z\in
Z_2\cap C\cap \Delta_n^{\circ}(x)\}$ est un ensemble de partitions
compatibles au sens de Burton et Keane, c'est-à-dire que, si
$P=\{P_1,P_2,P_3\}$ et $Q=\{Q_1,Q_2,Q_3\}$ appartienent à $\mathcal
P$, quitte à permuter les indices, on a $Q_2\cup Q_3\subset P_1$.
D'après le lemme de Burton et Keane dans \cite{BK}, on a
$\card\mathcal P\leq \card(\partial \Delta_n(x)\cap C)-2$ et notre
résultat en découle.\end{demo}

De ce lemme, on déduit immédiatement le

\begin{Cor}\label{Z2petit} Pour tout $x$ dans $X$ et pour tout entier naturel $n$, on a
$\card(Z_2\cap\Delta_n(x))\leq 5n+1$.\end{Cor}

Nous utilisons à présent ce lemme pour la

\begin{demo}[Démonstration de la proposition \ref{uniqueantecedent}] D'après le théorème de Birkhoff, il existe un point $x$ de $X$ tel qu'on ait
$$\frac{\card(Z_2\cap\Delta_n(x))}{\card(\Delta_n(x))}\td{n}{\infty}\mu(Z_2).$$
Or, d'après le corollaire \ref{Z2petit}, cette quantité tend vers $0$.\end{demo}

Finalement, nous utiliserons les résultats de cette section sous la forme du

\begin{Cor} \label{uniqueantecedent1}
Soit $\mu$ une probabilité borélienne sans atomes, invariante et ergodique par l'action de $\mathbb Z^2$ sur $X$. Alors, pour $\mu$-presque tout $x$ dans $X$,
pour toute $\sigma$-composante $C$ de la $\mathbb Z^2$-orbite de $x$, $C\cap Z$ contient au plus un point de la forme $T^kS^{-k}x$ pour $k$ dans $\mathbb Z$.\end{Cor}

\begin{demo}  Soit $U$ l'ensemble des éléments de $X$ qui ne satisfont pas la conclusion du corollaire et soit $x$ dans $U$. Alors, il existe des
entiers relatifs distincts $k$ et $l$ tels que les points $y=T^kS^{-k}x$ et $z=T^lS^{-l}x$ soient dans $Z$ et appartiennent à la même $\sigma$-composante.
Soit $q$ le plus petit entier tel que $\sigma^q(y)=\sigma^q(z)$. Alors, on a $\sigma^q(y)\in Z_2$. En d'autres termes, nous venons de montrer qu'on a $U\subset\bigcup_{i,j\in\mathbb Z}T^iS^jZ_2$. Il vient, d'après la proposition \ref{uniqueantecedent}, $\mu(U)=0$.
\end{demo}

\section{Les $\sigma$-composantes}
\label{compo}

Dans cette section, nous établissons un certain nombre de propriétés topologiques des $\sigma$-composantes.

Soient $x$ un point de $Y$ et $C$ sa $\sigma$-composante. Nous dirons que $x$ est un point extrémal de $C$ si les autres points de $C$ sont de la forme $T^kS^lx$ avec
$k+l>0$ ou $k+l=0$ et $l\geq 0$. Une $\sigma$-composante contient au plus un point extrémal. On note $E$ l'ensemble des éléments de $Y$ qui sont des points extrémaux de leur $\sigma$-composante. Nous avons la

\begin{Prop}\label{nonextremal} Soit $\mu$ une probabilité borélienne invariante et ergodique par l'action de $\mathbb Z^2$ sur $X$. On a $\mu(E)=0$.
\end{Prop}

Nous conservons les notations de la section \ref{antecedents}. La démonstration de la proposition \ref{nonextremal} découle encore du

\begin{Lem}\label{Epetit} Pour tout $x$ dans $X$ et pour tout entier naturel $n$, on a
$\card(E\cap\Delta_n(x))\leq (n+1)$.\end{Lem}

\begin{figure}\begin{center}\input{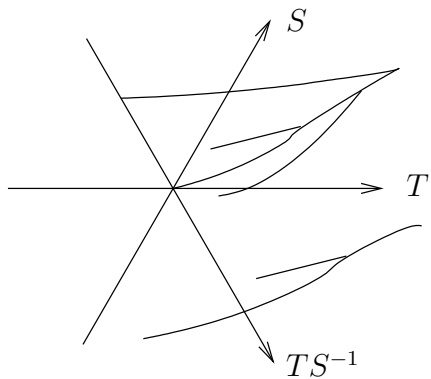}\caption{Un point extrémal}\label{extremal}\end{center}\end{figure}

\begin{demo} Soit $y$ un point de $E\cap\Delta_n(x)$. Alors, la $(\sigma,x,n)$-composante de $y$ contient un unique point de l'ensemble $\{T^kS^{n-k}x|0\leq k\leq n\}$. Ceci définit donc une application injective de $E\cap\Delta_n(x)$ dans un ensemble de cardinal $n+1$, d'où le résultat.\end{demo}

Soient $x$ dans $Y$ et $C$ sa $\sigma$-composante. Notons
\begin{align*}K(x)&=\{k\in\mathbb Z|T^kS^{-k}x\in C\}\\
L(x)&=\{k\in\mathbb Z|T^kS^{-k}x\in Y, T^kS^{-k}x\notin C\}\end{align*}
Alors, si $h\leq k\leq l$ sont des entiers relatifs tels que $h$ et $l$ appartiennent à $K(x)$ et que $T^kS^{-k}x$ appartienne à $Y$, on a $k\in K(x)$.

Notons $F$ l'ensemble des points $x$ de $Y$ pour lesquels la $\sigma$-composante de $x$ contient tous les points de la
$\mathbb Z^2$-orbite de $x$ qui appartiennent à $Y$
 et $G$ l'ensemble des $x$ dans $Y$ tels que $K(x)$ soit fini.
Nous avons la

\begin{Prop} \label{multicompo}
Soit $\mu$ une probabilité borélienne sans atomes, invariante et ergodique par l'action de $\mathbb Z^2$ sur $X$. On suppose qu'on a $\mu(F)=0$.
Alors, on a $\mu(G)=1$. En particulier, pour $\mu$-presque tout $x$ dans $X$, la $\mathbb Z^2$-orbite de $x$ contient une infinité de composantes.\end{Prop}

\begin{demo} Pour démontrer la proposition, vue notre remarque ci-dessus, il suffit d'établir que, pour $\mu$-presque tout $x$ dans $Y$,
l'ensemble $L(x)$ contient à la fois des éléments positifs et des éléments négatifs.

Commen\c cons par montrer que $L(x)$ n'est pas vide.
En effet, comme $\mu(F)=0$, pour $\mu$-presque tout $x$, il existe des entiers $k$ et $l$ tels que la $\sigma$-composante de $y=T^kS^lx$ ne contienne pas $x$.
Si $k+l\leq 0$, l'élément $\sigma^{-k-l}(y)$ est de la forme $T^hS^{-h}x$,
pour un certain entier $h$, et on a terminé puisqu'alors $h$ appartient à $L(x)$. Si $k+l>0$,
l'élément $z=\sigma^{k+l}(x)$ est de la forme $T^iS^jx$ avec $i+j=k+l$. Supposons par exemple qu'on a $i>k$.
La situation dans la $\mathbb Z^2$-orbite de $x$ des différents points entrant en jeu est alors représentée par la figure \ref{nouvellecompo}.

\begin{figure}\begin{center}\input{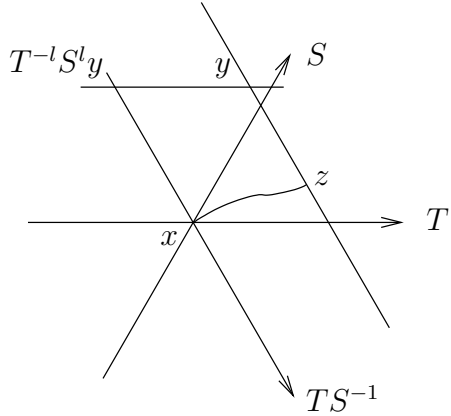}\caption{Démonstration de la proposition \ref{multicompo}} \label{nouvellecompo}\end{center}\end{figure}

Pour tout entier relatif $h$, si
$T^hS^{-h}y=T^{k+h}S^{l-h}x$ appartient à la $\sigma$-composante de $x$, on a $h>0$. Par conséquent, pour tout $m$, si $T^mS^{-m}x$ appartient à cette $\sigma$-composante, on a $m>-l$. Or, d'après le théorème de récurrence de Poincaré, il existe une infinité d'entiers positifs $m$ tels que $T^{-m}S^{m}x$ appartienne à $Y$. Pour $m$ suffisament grand, la $\sigma$-composante de $T^{-m}S^{m}x$ est donc différente de celle de $x$ : l'ensemble $L(x)$ n'est donc pas vide.

Considérons à présent, par exemple, l'ensemble $H$ des $x$ dans $Y$
tels que $L(x)$ ne contienne que des entiers positifs. Alors, d'après le raisonnement ci-dessus,
si $H'$ désigne l'ensemble des $x$ dans $H$ tels que, pour tout $k>0$, si $T^kS^{-k}x$ est dans $Y$,
$T^kS^{-k}x$ n'appartient pas à la $\sigma$-composante de $x$,
on a $\mu\left(H-\bigcup_{k\geq 0}T^{-k}S^kH'\right)=0$. Or, pour tout $x$ dans $H'$, pour tout $k>0$, on a  $T^kS^{-k}x\not\in H'$ et, donc, d'après le théorème de récurrence de Poincaré, $\mu(H')=0$. Il vient $\mu(H)=0$ : en d'autres termes, pour $\mu$-presque tout $x$, $L(x)$ contient des éléments négatifs.
\end{demo}

Des propositions \ref{Epetit} et \ref{multicompo}, on déduit le

\begin{Cor} \label{Zgros}
Soit $\mu$ une probabilité borélienne sans atomes, invariante et ergodique par l'action de $\mathbb Z^2$ sur $X$. On suppose qu'on a $\mu(F)=0$. Alors, pour $\mu$-presque tout $x$ dans $X$, la $\sigma$-composante de $x$ contient un point de $Z$.\end{Cor}

\begin{demo} Considérons l'ensemble $V$ des $x$ dans $X$ dont la $\sigma$-com\-po\-san\-te ne contient pas de points de $Z$. D'après la proposition \ref{multicompo}, pour $\mu$-presque tout $x$ dans $V$, la $\sigma$-composante $C$ de $x$ ne contient qu'un nombre fini de points de la forme $T^kS^{-k}x$ pour $k$ dans $\mathbb Z$. Mais alors, comme chacun de ces points ne possède qu'un nombre fini de $\sigma$-antécédents, $C$ ne contient qu'un nombre fini de points de la forme $T^kS^lx$ avec $k+l\leq 0$. En particulier, $C$ contient un point extrémal. En d'autres termes, on a $\mu\left(V-\bigcup_{k,l\in\mathbb Z}T^{k}S^lE\right)=0$. D'après la proposition \ref{Epetit}, on a donc $\mu(V)=0$.
\end{demo}

\section{Démonstration du théorème}
\label{demothm}

Nous séparons les deux cas qui apparaissent dans le théorème.

\begin{demo}[Démonstration du théorème dans le cas où l'on a $\mu(Z)=0$] Il s'agit de montrer qu'on a, avec les notations de la section \ref{compo}, $\mu(F)=1$.
Dans le cas contraire, par ergodicité, on a $\mu(F)=0$ et, donc, d'après le corollaire \ref{Zgros}, $\mu\left(\bigcup_{k,l\in\mathbb Z}T^{k}S^lZ\right)=1$. Il vient bien $\mu(F)=1$.\end{demo}

\begin{demo}[Démonstration du théorème dans le cas où l'on a $\mu(Z)>0$] Commen\c cons par montrer qu'on a $\mu(F)=0$. Supposons, au contraire, qu'on a $\mu(F)=1$. Alors, d'après le théorème de récurrence de Poincaré, pour $\mu$-presque tout $x$ dans $Z\cap F$, il existe un entier $k\neq 0$ tel que $y=T^kS^{-k}x$ appartienne à $Z$.
Mais alors, $x$ et $y$ appartiennent à la même $\sigma$-composante. D'après le corollaire \ref{uniqueantecedent1}, on a donc $\mu(Z\cap F)=0$, ce qui est contradictoire.

D'après la proposition \ref{multicompo} et le corollaire \ref{Zgros}, pour $\mu$-presque tout $x$ dans $X$, la $\mathbb Z^2$-orbite de $x$ contient une infinité de $\sigma$-composantes et chacune d'entre elles contient un point de $Z$. En particulier, chacune d'entre elles contient un point de $Z$ qui est de la forme $T^kS^{-k}x$ pour un $k$ dans $\mathbb Z$. D'après
le corollaire \ref{uniqueantecedent1}, pour $\mu$-presque tout $x$, pour toute $\sigma$-composante de la $\mathbb Z^2$-orbite de $x$, ce point est unique. Enfin, d'après la proposition \ref{uniqueantecedent}, pour $\mu$-presque tout $x$ dans $Z$, $\sigma^{-1}(x)$ contient un unique point de $Z$.
\end{demo}

\section{Une mesure du type arbre}
\label{exarbre}

On considère dorénavant $(\mathbb Z/2\mathbb Z)^{\mathbb Z^2}$ comme un groupe compact pour la loi produit et $X$ comme un sous-groupe fermé de $(\mathbb Z/2\mathbb Z)^{\mathbb Z^2}$. Alors, $\mathbb Z^2$ agit par automorphismes de groupe sur $X$. En particulier, cette action préserve la mesure de Haar $\mu_0$ de $X$. Nous pouvons décrire plus précisément cette mesure. Notons $P\subset\mathbb Z^2$ la réunion des ensembles $P_-=\{(k,-k)|k\in\mathbb Z\}$ et $P_+=\{(l,0)|l\in\mathbb N^*\}$. On a immédiatement le

\begin{Lem}\label{Haar} L'application $X\rightarrow (\mathbb Z/2\mathbb Z)^{P},x\mapsto (x_{k,l})_{(k,l)\in P}$ est un isomorphisme de groupes compacts.
\end{Lem}

Par conséquent, en termes probabilistes, pour tirer au hasard un élément $x$ de $X$ suivant la loi $\mu_0$, on tire au hasard de manière indépendante suivant la loi de Bernoulli de paramètre $\frac{1}{2}$ les coordonnées $(x_{k,l})_{(k,l)\in P}$ de $x$ et on complète de proche en proche grâce à l'équation $x_{k,l}+x_{k+1,l}+x_{k,l+1}=0$.

Un argument élémentaire de transformation de Fourier dans le groupe abélien compact $X$ permet de dé\-mon\-trer le

\begin{Lem}\label{melange} La mesure $\mu_0$ est globalement mélangeante pour l'action de $\mathbb Z^2$. Plus précidément, pour tous boréliens $A$ et $B$ de $X$, on a
$$\mu_0(A\cap T^kS^lB)\td{(k,l)}{\infty}\mu_0(A)\mu_0(B).$$\end{Lem}

En particulier, la mesure $\mu_0$ est ergodique.

Nous allons démontrer la

\begin{Prop}\label{Haararbre} La mesure $\mu_0$ est du type arbre. \end{Prop}

La démonstration repose sur l'argument suivant de prolongement :

\begin{Lem}\label{prolongement} Soit $y=(y_{k,l})_{k+l\leq 0}$ une famille d'éléments de $\mathbb Z/2\mathbb Z$ telle que, pour tous entiers relatifs $k$ et $l$ avec $k+l+1\leq 0$, on ait
$y_{k,l}+y_{k+1,l}+y_{k,l+1}=0$. On suppose qu'on a $y_{0,0}=1$ et, pour un certain $h\geq 0$, $y_{h,-h}=1$. Alors, il existe $x$ dans $X$ tel que, pour tous $k,l$ avec
$k+l\leq 0$, $x_{k,l}=y_{k,l}$ et que $x$ et $T^hS^{-h}x$ soient dans la même $\sigma$-composante.\end{Lem}

En d'autres termes, étant donné le passé d'un élément de $X$, on peut, en changeant son futur, forcer deux éléments de son orbite à se rejoindre par $\sigma$.

\begin{demo} La démonstration se fait par récurrence sur $h\geq 0$. Si $h=0$, le résultat est évident.

Supposons donc $h\geq 1$ et le résultat établi pour $h-1$. Et supposons, par l'absurde, que $y$ ne possède pas de prolongement convenable.
Considérons une famille $z=(z_{k,l})_{k+l\leq 1}$ qui prolonge $y$ et telle que, pour tous $k$ et $l$ avec
$k+l\leq 0$, on ait $z_{k,l}+z_{k+1,l}+z_{k,l+1}=0$. Alors $z$ est complètement déterminée par la donnée de $z_{1,0}$ ; l'unique autre prolongement de $y$ à l'ensemble
$\{k+l\leq 1\}$ qui vérifie la même condition a pour valeur $z_{k,1-k}+1$ en $(k,1-k)$ pour tout $k$ dans $\mathbb Z$. Deux cas sont donc possibles :
soit le prolongement de $y$ valant $1$ en $(1,0)$ a pour valeur $1$ en $(h,1-h)$ et, par récurrence, on pourrait alors prolonger cette famille en une famille
$x$ qui convient, ce qu'on a supposé impossible ;
soit le prolongement de $y$ valant $1$ en $(0,1)$ vaut $1$ en $(h,1-h)$. En réitérant ce raisonnement, on peut construire, pour tout $1\leq i\leq h$, une famille
$z^{(i)}=\left(z^{(i)}_{k+l}\right)_{k+l\leq i}$ de prolongements successifs de $y$ telle que, pour tout $i$ et pour tous $k,l$ avec $k+l\leq i-1$, on ait
$z^{(i)}_{k,l}+z^{(i)}_{k+1,l}+z^{(i)}_{k,l+1}=0$ et  $z^{(i)}_{0,i}=z^{(i)}_{h,i-h}=1$. Mais alors, pour $i=h$, on a $z^{(h)}_{h,0}=1$ tandis que, pour tout $0\leq k\leq h$,
$z^{(h)}_{0,k}=1$, ce qui implique  que $z^{(h)}_{h,0}=0$. On arrive à une contradiction, le résultat en découle.
\end{demo}

\begin{demo}[Démonstration de la proposition \ref{Haararbre}]
Supposons au contraire que la mesure $\mu_0$ soit du type ruban et
donnons-nous un ensemble $B$ comme dans le théorème. Pour $x$ dans
$X$ et $y$  dans $(\mathbb Z/2\mathbb Z)^{P_+}$, nous noterons
$[x,y]$ l'unique élément $z$ de $X$ tel $z_{k,l}=x_{k,l}$ pour
$k+l\leq 0$ et que $z_{k,0}=y_{k,0}$ pour $k>0$. Alors, comme,
d'après le lemme \ref{Haar}, $\mu_0$ s'identifie à la mesure produit
des mesures de Haar $\mu_-$ de $(\mathbb Z/2\mathbb Z)^{P_-}$ et
$\mu_+$ de  $(\mathbb Z/2\mathbb Z)^{P_+}$, on peut trouver un
élément $x$ de $B\cap Z$ tel que, pour $\mu_+$-presque tout $y$ dans
$(\mathbb Z/2\mathbb Z)^{P_+}$, l'élément $[x,y]$ appartienne à $B$.
En particulier, comme $x$ est dans $B$, il existe un entier $h>0$
tel que $T^hS^{-h}x$ appartienne à $Z$. D'après le lemme
\ref{prolongement}, il existe un entier $l>0$ et des éléments
$t_1,\ldots,t_l$ de $\mathbb Z/2\mathbb Z$ tels que, pour tout $y$
dans $(\mathbb Z/2\mathbb Z)^{P_+}$, si
$y_{1,0}=t_1,\ldots,y_{l,0}=t_l$, les points $[x,y]$ et
$T^hS^{-h}[x,y]$ appartiennent à la même $\sigma$-composante. Mais
alors, si $[x,y]$ appartient à $B$, ceci contredit le fait que la
$\sigma$-composante de $[x,y]$ ne contient qu'un seul point de
$Z\cap\{T^kS^{-k}[x,y]|k\in\mathbb Z\}$.
\end{demo}

\section{Une mesure du type ruban}
\label{exruban}

Dans cette section, nous allons construire une mesure du type ruban.

Commen\c cons par considérer l'ensemble $X_1=\{x\in X|\forall k,l\in\mathbb Z\quad x_{2k,2l}=0\}$ : c'est un sous-groupe fermé de $X$ qui est stable par l'action de $(2\mathbb Z)^2$. Pour $x$ dans $X$, et pour $k$ et $l$ dans $\mathbb Z$, posons $\theta_{2k,2l}(x)=0$ et $\theta_{2k+1,2l}(x)=\theta_{2k,2l+1}(x)=\theta_{2k+1,2l+1}(x)=x_{k,l}$. On vérifie aisément que l'application $\theta:x\mapsto (\theta_{k,l}(x))_{(k,l)\in\mathbb Z^2}$ définit un isomorphisme de $X$ sur $X_1$ et que, pour tous entiers $k$ et $l$, on a
$\theta T^{2k}S^{2l}=T^kS^l\theta$. En particulier, d'après le lemme \ref{melange}, la mesure $\theta_*\mu_0$ est invariante et ergodique pour l'action de $(2\mathbb Z)^2$ sur $X_1$. On pose $\mu_1=\frac{1}{4}(\theta_*\mu_0+T_*\theta_*\mu_0+S_*\theta_*\mu_0+(TS)_*\theta_*\mu_0)$ : la mesure $\mu_1$ est invariante et ergodique pour l'action de $\mathbb Z^2$ sur $X$ et son support est l'ensemble $X_1\cup TX_1\cup SX_1 \cup TSX_1$. Enfin, comme $X_1\cap TX_1=X_1\cap SX_1=X_1\cap TSX_1=\{0\}$,
il existe une application $\pi$ du support de $\mu_1$ privé de $0$ dans $(\mathbb Z/2\mathbb Z)^2$ telle que $\pi(X_1)=(0,0)$, $\pi(TX_1)=(1,0)$, $\pi(SX_1)=(0,1)$ et $\pi(TSX_1)=(1,1)$. En particulier, l'application $\pi$ entrelace l'action de $\mathbb Z^2$ sur le support de $X_1$ et l'action naturelle par translations de $\mathbb Z^2$ sur
$(\mathbb Z/2\mathbb Z)^2$.

Par ailleurs, pour $k$ et $l$ dans $\mathbb Z$, posons $u_{k,l}=0$ si $k-l$ est congru à $0$ modulo $3$ et $u_{k,l}=1$ sinon. Alors la famille $u=(u_{k,l})$ appartient à $X$ et $u$ est un point périodique pour l'action de $\mathbb Z^2$ sur $X$ : son stabilisateur est précisément l'ensemble des $(k,l)$ dans $\mathbb Z^2$ tels que $k-l$ soit congru à $0$ modulo $3$. Pour $v$ dans $\{u,Tu,Su\}$, on pose $\kappa(v)=0$ si $v=u$, $\kappa(v)=1$ si $v=Tu$ et $\kappa(v)=2$ si $v=Su$. Alors, l'application $\kappa:\{u,Tu,Su\}\rightarrow \mathbb Z/3\mathbb Z$ entrelace l'action de $\mathbb Z^2$ sur la $\mathbb Z^2$-orbite de $u$ et l'action de
$\mathbb Z^2$ sur $\mathbb Z/3\mathbb Z$ pour laquelle, pour tout $(k,l)$ dans $\mathbb Z^2$, $(k,l)$ agit sur $\mathbb Z/3\mathbb Z$ par la translation par $k-l$.
On note $\nu$ la mesure de probabilité invariante portée par la $\mathbb Z^2$-orbite de $u$, c'est-à-dire la mesure $\frac{1}{3}(\delta_u+\delta_{Tu}+\delta_{Su})$, où $\delta$ désigne une masse de Dirac.

Enfin, on pose $\mu=\mu_1*\nu$, le produit de convolution des deux mesures de probabilité $\mu_1$ et $\nu$ sur le groupe compact $X$, c'est-à-dire l'image par l'application somme de la mesure $\mu_1\otimes\nu$ sur le produit cartésien $X\times X$. Comme $\mathbb Z^2$ agit sur $X$ par automorphismes de groupe, la mesure $\mu$ est encore invariante par cette action.

\begin{Lem} \label{structureproduit}
L'application produit $(X\times X,\mu_1\otimes\nu)\rightarrow (X,\mu)$ est un isomorphisme mesuré. La mesure $\mu$ est ergodique pour l'action de $\mathbb Z^2$.
\end{Lem}

\begin{demo} Commen\c cons par montrer la première assertion. Le
groupe diédral d'ordre $6$ agit de manière naturelle sur notre
situation. On vérifie que, par symétrie, il suffit de montrer qu'on
a $\mu((u+X_1)\cap(u+TX_1))
=\mu((u+X_1)\cap(Tu+X_1))=\mu((u+X_1)\cap(Tu+TX_1))=0$. Comme
$X_1\cap TX_1=\{0\}$, on a $\mu((u+X_1)\cap(u+TX_1))=0$. Comme
$(Tu-u)_{0,0}=1$, on a $(u+X_1)\cap(Tu+X_1)=\emptyset$. Soit enfin
$x$ dans $(u+X_1)\cap(Tu+TX_1)$. Comme $x$ appartient à $u+X_1$, on
vérifie qu'on a, pour tout $l$, pour tout $k$ congru à $2$ ou $4$
modulo $6$, $x_{k,2l}=1$. De même, comme $x$ appartient à $Tu+TX_1$,
on vérifie qu'on a, pour tout $l$, pour tout $k$ congru à $1$ ou $3$
modulo $6$, $x_{k,2l}=1$. Alors, on a en particulier,
$x_{1,0}=x_{2,0}=x_{3,0}=1$, donc $x_{1,1}=x_{2,1}=0$ et
$x_{1,2}=0$, ce qui est contradictoire. Il vient
$(u+X_1)\cap(Tu+TX_1)=\emptyset$.

Le système $(X,\mathbb Z^2,\mu)$ est donc exactement isomorphe au système produit $(X,\mathbb Z^2,\mu_1)\times (X,\mathbb Z^2,\nu)$. D'après le lemme \ref{melange}, le spectre discret de $(X,\mathbb Z^2,\mu_1)$ est égal à $(\frac{1}{2}\mathbb Z/\mathbb Z)^2\subset \mathbb T^2=(\mathbb R/\mathbb Z)^2$. Par ailleurs, par un calcul direct, le spectre discret de $(X,\mathbb Z^2,\nu)$ est égal à $\{(1,1),(\frac{1}{3},\frac{2}{3}),(\frac{2}{3},\frac{1}{3})\}\subset \mathbb T^2$. Comme ces deux sous-groupes ne s'intersectent qu'en $1$, par un argument classique de théorie ergodique, le système produit est ergodique.
\end{demo}

Nous avons la

\begin{Prop} La mesure $\mu$ est du type ruban.\end{Prop}

\begin{demo} D'après le lemme \ref{structureproduit},
il existe une unique application définie $\mu$-presque partout sur $X$, $x\mapsto(y(x),v(x)),X\rightarrow X\times \{u,Tu,Su\}$ telle  que $y_*\mu=\mu_1$, que $v_*\mu=\nu$ et que, pour $\mu$-presque tout point $x$ de $X$, on ait $x=y(x)+v(x)$. On pose alors, pour $\mu$-presque tout $x$ dans $X$,
$\varpi(x)=(\pi(y(x)),\kappa(v(x)))\in(\mathbb Z/2\mathbb Z)^2\times\mathbb Z/3\mathbb Z$.
On notera $00$, $10$, $01$ et $11$ les quatre éléments de  $(\mathbb Z/2\mathbb Z)^2$.

Pour $\mu$-presque tout $x$ dans $Y$, si $\varpi(x)=(a,b)$, on a
$\varpi(\sigma(x))=(a+10,b+1)$ ou $\varpi(\sigma(x))=(a+01,b+2)$,
suivant que $\sigma(x)=Tx$ ou que $\sigma(x)=Sx$. Il existe des
$(a,b)$ dans $(\mathbb Z/2\mathbb Z)^2\times\mathbb Z/3\mathbb Z$
tels que, pour $\mu$-presque tout $x$ dans $Y$, si
$\varpi(x)=(a,b)$, on a nécessairement $\sigma(x)=Tx$ ou
$\sigma(x)=Sx$. Par exemple, si $\varpi(x)=(01,0)$, on a, d'après la
définition de $\mu_1$ et de $\nu$,  $y_{1,0}(x)=y_{0,0}(x)=1$ et
$y_{0,1}(x)=0$, tandis que
$v_{1,0}(x)=u_{1,0}=1=u_{0,1}=v_{0,1}(x)$, si bien que $x_{1,0}=0$
et $x_{0,1}=1$ et, donc, $\sigma(x)=Sx$. De même, si
$\varpi(x)=(10,0)$, on a $\sigma(x)=Tx$. Sur la figure \ref{marche},
on a représenté l'ensemble de ces contraintes : pour chaque $(a,b)$
dans $(\mathbb Z/2\mathbb Z)^2\times\mathbb Z/3\mathbb Z$, on a
indiqué, par une ou deux flèches, les valeurs que peut prendre
$\varpi(\sigma(x))$ quand $\varpi(x)$ est donné. \`A la base de
chaque flèche, on a fait figurer la lettre $T$ ou $S$ suivant que la
valeur atteinte l'est quand $\sigma(x)=Tx$ ou quand $\sigma(x)=Sx$.

\begin{figure}\begin{center}\input{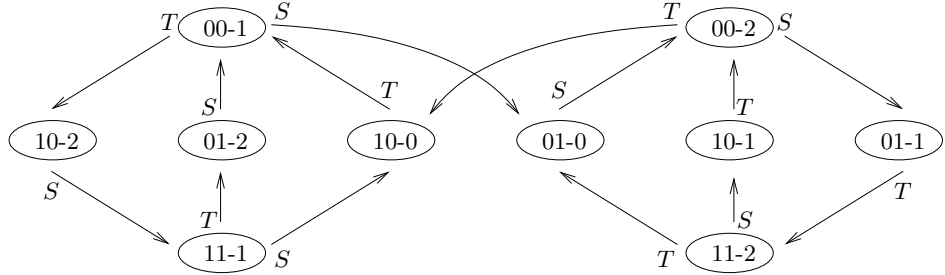}\caption{Contraintes sur les valeurs successives des $\varpi\circ\sigma^n$, $n\in\mathbb N$}
\label{marche}\end{center}\end{figure}

Soient $\varphi:(\mathbb Z/2\mathbb Z)^2\times\mathbb Z/3\mathbb
Z\rightarrow \mathbb Z$ avec
\begin{align*}\varphi(a,b)&=-2\mbox{ pour }(a,b)\in\{(10,1),(01,1)\}\\
\varphi(a,b)&=-1\mbox{ pour }(a,b)\in\{(00,2),(11,2)\}\\
\varphi(a,b)&=0\mbox{ pour }(a,b)\in\{(10,0),(01,0)\}\\
\varphi(a,b)&=1\mbox{ pour }(a,b)\in\{(00,1),(11,1)\}\\
\varphi(a,b)&=2\mbox{ pour }(a,b)\in\{(10,2),(01,2)\}
\end{align*}
et $\psi=\varphi\circ\varpi$.
On vérifie que, d'après la figure \ref{marche},
pour $\mu$-presque tout $x$ dans $Y$, on a $\sigma(x)=Tx$ si $\psi(\sigma(x))-\psi(x)=1$ et
$\sigma(x)=Sx$ si $\psi(\sigma(x))-\psi(x)=-1$. En d'autres termes, pour $\mu$-presque tout $x$ dans $Y$, on a
$$\sigma(x)=T^{\frac{1}{2}(1+\psi(\sigma(x))-\psi(x))}S^{\frac{1}{2}(1-\psi(\sigma(x))+\psi(x))}x.$$
Par récurrence, il vient, pour tout $n$ dans $\mathbb N$,
$$\sigma^n(x)=T^{\frac{1}{2}(n+\psi(\sigma^n(x))-\psi(x))}S^{\frac{1}{2}(n-\psi(\sigma^n(x))+\psi(x))}x$$
et, donc, comme $\psi$ prend toute ses valeurs dans $[-2,2]$, $\sigma^n(x)$ est de la forme $T^kS^lx$ avec $k+l=n$ et $\abs{k-l}\leq 4$. Par conséquent, pour $\mu$-presque tout $x$ dans $Y$, pour tout $h$ dans $\mathbb Z$ avec $\abs{h}>8$, si $T^hS^{-h}x$ appartient à $Y$, les $\sigma$-composantes de $x$ et de  $T^hS^{-h}x$ sont disjointes. La mesure $\mu$ est donc du type ruban.
\end{demo}

\noindent Jean-François Quint\\
LAGA\\
Université Paris 13\\
99, avenue Jean-Baptiste Clément\\
93430 Villetaneuse\\
France\\
quint@math.univ-paris13.fr


\begin{thebibliography}{1}

\bibitem{BK} R. M. Burton, M. Keane, Density and uniqueness in percolation, {\em Communications in mathematical physics} {\bf 121} (1989), 501-505.

\bibitem{Fu} H. Fürstenberg, Disjointness in ergodic theory, minimal sets, and a problem in Diophantine approximation,
{\em  Mathematical systems theory}  {\bf 1} (1967), 1-49.

\bibitem{HMM} B. Host, A. Maass, S. Martínez, Uniform Bernoulli measure in dynamics of permutative cellular automata with algebraic local rules, {\em Discrete and continuous dynamical systems} {\bf 9} (2003), 1423-1446.

\bibitem{Sab} M. Sablik, Measure rigidity for algebraic bipermutative cellular automata, preprint.

\bibitem{Silb} S. Silberger, Subshifts of the three dot system, {\em Ergodic theory and dynamical systems} {\bf 25}  (2005), 1673-1687.

\end{thebibliography}
\end{document}